\documentclass[draft]{jotart}
\usepackage{amsfonts}
\usepackage{amsmath}
\usepackage{amssymb}
\usepackage{mathrsfs}
\theoremstyle{proclaim}
\newtheorem{theorem}{Theorem}[section]
\newtheorem{lemma}[theorem]{Lemma}
\newtheorem{corollary}[theorem]{Corollary}

\theoremstyle{statement}

\newtheorem{example}[theorem]{Example}
\theoremstyle{fancyproclaim}

\numberwithin{equation}{section}

\begin{document}
\issueinfo{00}{0}{0000}
\commby{Editor}
\pagespan{101}{110}
\date{Month dd, yyyy}
\revision{Month dd, yyyy}

\title[On the Strongly Irreducible Decompositions of Operators]
{On The Uniqueness of The Strongly Irreducible Decompositions
of Operators up to Similarity} 
\author[Rui Shi]{Rui Shi}
\address{RUI SHI, School of Mathematical Sciences, Dalian University
of Technology, Dalian, 116024, China}
\email{littlestoneshg@gmail.com}

\begin{abstract}
We give a generalization of the Jordan canonical form theorem for a
class of bounded linear operators on complex separable Hilbert
spaces in terms of direct integrals. Precisely, we study the
uniqueness of strongly irreducible decompositions of the operators
on the Hilbert spaces up to similarity.
\end{abstract}


\maketitle

\section{INTRODUCTION}

Throughout this article, all Hilbert spaces discussed are
{\textit{complex and separable}}. Denote by
$\mathscr{L}(\mathscr{H})$ the set of bounded linear operators on a
Hilbert space $\mathscr{H}$. An {\textit{idempotent}} $P$ is an
operator in $\mathscr{L}(\mathscr{H})$ satisfying $P^{2}=P$. A
{\textit{projection}} $Q$ is an idempotent such that $\textrm{ker}
Q=(\textrm{ran} Q)^{\perp}$ (See \cite{Conway_1}). An operator $A$
in $\mathscr {L}(\mathscr {H})$ is said to be {\textit{irreducible}}
if its commutant $\{A\}'\triangleq\{B\in\mathscr {L}(\mathscr
{H}):AB=BA\}$ contains no projections other than $0$ and the
identity operator $I$ on $\mathscr{H}$, introduced by P. Halmos in
\cite{Halmos}. (The separability assumption is necessary because on
a non-separable Hilbert space every operator is reducible.) An
operator $A$ in $\mathscr {L}(\mathscr {H})$ is said to be
{\textit{strongly irreducible}} if $XAX^{-1}$ is irreducible for
every invertible operator $X$ in $\mathscr {L}(\mathscr {H})$
\cite{Gilfeather}. This shows that the commutant of a strongly
irreducible operator contains no idempotents other than $0$ and $I$.
Strong irreducibility stays invariant up to similar equivalence
while irreducibility is only an invariant up to unitary equivalence.
An idempotent $P$ in $\{A\}^{\prime}$ is said to be
{\textit{minimal}} if every idempotent $Q$ in
$\{A\}^{\prime}\cap\{P\}^{\prime}$ satisfies $QP=P$ or $QP=0$. For a
minimal idempotent $P$ in $\{A\}^{\prime}$, it can be observed that
the restriction $A|_{\textrm{ran} P}$ is strongly irreducible on
$\textrm{ran} P$. An operator $A$ in $\mathscr {L}(\mathscr {H})$ is
said to have a finite {\textit{strongly irreducible decomposition}}
if there exist finitely many minimal idempotents
$\{P_{i}\}^{n}_{i=1}$ in $\{A\}^{\prime}$ such that
$\sum^{n}_{i=1}P_{i}=I$ and $P_{i}P_{j}=P_{j}P_{i}=0$ for $1\leq
i\neq j\leq n$. By the above observation, an operator $A$ in
$\mathscr {L}(\mathscr {H})$ having a finite strongly irreducible
decomposition can be expressed as a direct sum of finitely many
strongly irreducible operators.

On finite dimensional Hilbert spaces, every strongly irreducible
operator is similar to a Jordan block. In \cite{Herrero_1}, D. A.
Herrero and C. Jiang proved that for every operator $T$ in $\mathscr
{L}(\mathscr {H})$, there exists a sequence
$\{T_{n}\}^{\infty}_{n=1}$ in $\mathscr {L}(\mathscr {H})$ such that
$\lim_{n\rightarrow\infty}\|T-T_{n}\|=0$, where every operator
$T_{n}$ is similar to a direct sum of finitely many strongly
irreducible operators. Y. Cao, J. Fang and C. Jiang \cite{Jiang_2}
studied the uniqueness of finite strongly irreducible decompositions
of operators in $\mathscr {L}(\mathscr {H})$ up to similar
equivalence by the $K_{0}$ groups of Banach algebras. For more work
around this subject, the reader is referred to \cite{Watatani_1,
Watatani_2, Watatani_3, Jiang_3, Jiang_1, Jiang_4, Jiang_5, Jiang_6,
Jiang Z}. Inspired by the ideas and results in \cite{Jiang_2}, we
study operators in $\mathscr {L}(\mathscr {H})$ which may have no
finite strongly irreducible decompositions. In particular, there are
many operators in $\mathscr {L}(\mathscr {H})$ whose commutants
contain no minimal idempotents. To represent these operators, direct
sums of strongly irreducible operators need to be generalized to
direct integrals with some regular Borel measures. In \cite{Shi_1},
C. Jiang and the author of the present paper proved that an operator
$A$ in $\mathscr {L}(\mathscr {H})$ is similar to a direct integral
of strongly irreducible operators if and only if its commutant
$\{A\}^{\prime}$ contains a bounded maximal abelian set of
idempotents. A direct integral of strongly irreducible operators
means the integrand is strongly irreducible almost everywhere on the
domain of integration. For related concepts and results about direct
integrals and abelian von Neumann algebras, the reader is referred
to \cite{Azoff_2, Conway_1, Davidson, Rosenthal, Schwartz}.

Following the notation of \cite{Shi_1}, we generalize a definition
mentioned above. An operator $A$ in $\mathscr {L}(\mathscr {H})$ is
said to have a {\textit{strongly irreducible decomposition}} if its
commutant $\{A\}^{\prime}$ contains a bounded maximal abelian set of
idempotents. Furthermore, a strongly irreducible decomposition of
the operator $A$ is said to be {\textit{unique up to similarity}} if
for bounded maximal abelian sets of idempotents $\mathscr{P}$ and
$\mathscr{Q}$ in $\{A\}^{\prime}$, there is an invertible operator
$X$ in $\{A\}^{\prime}$ such that $X\mathscr{P}X^{-1}=\mathscr{Q}$.

As a corollary of the main theorems, a normal operator in $\mathscr
{L}(\mathscr {H})$ has unique strongly irreducible decomposition up
to similarity if and only if the multiplicity function $m_{_N}$ for
$N$ is finite a.~e.\ on $\sigma(N)$ with respect to the
scalar-valued spectral measure $\mu_{_N}$. By this, the tensor
product $I_{\mathscr{H}}\otimes N$ does not have unique strongly
irreducible decomposition up to similarity, if
${\textrm{dim}}\mathscr{H}=\infty$.

To simplify the statements of the main theorems, we need to
introduce the upper triangular representation for operator-valued
matrices. Assume $A$ in $\mathscr {L}(\mathscr {H})$ is a direct
integral of strongly irreducible operators in the form
$$A=\bigoplus^{\infty}_{n=1}\int_{\Lambda_n} A(\lambda)
d\mu(\lambda), \eqno (1)$$ with respect to a partitioned measure
space $\{\Lambda, \mu, \{\Lambda_{n}\}^{\infty}_{n=1}\}$, where
$\mu$ is a regular Borel measure on a compact set $\Lambda$ and
$\{\Lambda_{n}\}^{\infty}_{n=1}$ is a Borel partition of $\Lambda$,
and the equation $\mu(\Lambda_{n})=0$ holds for all but finitely
many $n$ in $\mathbb{N}$ ($0\notin\mathbb{N}$), and the dimension of
the fibre space $\mathscr{H}_{\lambda}$ (\cite{Abrahamse}, \S2) is
$n$ for almost every $\lambda$ in $\Lambda_{n}$.

By (\cite{Azoff_1}, Corollary 2), there is a unitary operator $U$
such that $$UAU^{*}=
\bigoplus\limits^{\infty}_{n=1}\int\limits_{\Lambda_n}
\left(\begin{array}{cccccc}
M_{\phi_{n}}&M_{\phi^{n}_{12}}&M_{\phi^{n}_{13}}&\cdots&M_{\phi^{n}_{1n}}\\
0&M_{\phi_{n}}&M_{\phi^{n}_{23}}&\cdots&M_{\phi^{n}_{2n}}\\
0&0&M_{\phi_{n}}&\cdots&M_{\phi^{n}_{3n}}\\
\vdots&\vdots&\vdots&\ddots&\vdots\\
0&0&0&\cdots&M_{\phi_{n}}\\
\end{array}\right)_{n\times n}(\lambda) d\mu_{n}(\lambda),\eqno (2)$$
where $\mu_{n}=\mu|_{\Lambda_{n}}$ for $1\leq n<\infty$ and
$\phi_{n},\phi^{n}_{ij}\in L^{\infty}(\mu_{n})$, and $M_{\phi_{n}}$,
$M_{\phi^{n}_{ij}}$ are multiplication operators. Denote by
$\nu_{n}=\mu_{n}\circ\phi^{-1}_{n}$ the scalar-valued spectral
measure for $M_{\phi_{n}}$. Let the set
$\{\Gamma_{nm}\}^{m=\infty}_{m=1}$ be the Borel partition of the
spectrum $\sigma(M_{\phi_{n}})$ with respect to the
$\nu_{n}$-measurable multiplicity function $m_{_{{\phi_{_n}}}}$ for
$M_{\phi_{n}}$ defined on $\sigma(M_{\phi_{n}})$ such that
$m_{_{{\phi_{_n}}}}(\lambda)=m$ for almost every $\lambda$ in
$\Gamma_{nm}$. Write $\nu_{nm}$ for $\nu_{n}|_{\Gamma_{nm}},1\leq
m\leq\infty$.

For a class of operators in $\mathscr {L}(\mathscr {H})$ having
unique strongly irreducible decompositions up to similarity, we give
a necessary and sufficient condition by $K$-theory for Banach
algebras. Precisely, we prove the following theorems.

\begin{theorem}
Assume that an operator $A$ in $\mathscr {L}(\mathscr {H})$ is
stated as in (1) and expressed as in (2) such that
\begin{enumerate}
\item the $\nu_{n}$-measurable multiplicity function
$m_{_{{\phi_{_n}}}}$ is simple and may take $\infty$ on the spectrum
$\sigma(M_{\phi_{n}})$ for every $n$ in $\mathbb{N}$ and
\item every superdiagonal entry as in (2) is invertible for
$n$ in $\{n\in\mathbb{N}:\mu(\Lambda_{n})>0\}$.
\end{enumerate}
Then the following statements are equivalent.
\begin{enumerate}
\item [(a)] The strongly irreducible decomposition of $A$ is
unique up to similarity.
\item [(b)] There exists a bounded $\mathbb{N}$-valued simple
function $r_{_A}$ on $\sigma(A)$ such that
\begin{enumerate}
\item [~] $V(\{A\}^{\prime})\cong\{f(\lambda)\in\mathbb{N}^{(r_{_A}(\lambda))}:f\
{\textrm{is\ Borel and bounded on }}\sigma(A)\}$ and
\item [~] $K_{0}(\{A\}^{\prime})\cong\{f(\lambda)\in\mathbb{Z}^{(r_{_A}(\lambda))}:f\
{\textrm{is\ Borel and bounded on }}\sigma(A)\}.$
\end{enumerate}
\end{enumerate}
\end{theorem}

\begin{theorem}
If an operator $A$ in $\mathscr {L}(\mathscr {H})$ is expressed as
in (1) and (2) such that the $\nu_{n}$-measurable multiplicity
function $m_{_{{\phi_{_n}}}}$ is simple and bounded on
$\sigma(M_{\phi_{n}})$ for every $n$ in $\mathbb{N}$, then there
exists a sequence of operators $\{A_{k}\}^{\infty}_{k=1}$ in
$\mathscr {L}(\mathscr {H})$ required as in Theorem 1.1 and having
unique strongly irreducible decompositions up to similarity such
that $\lim_{k\rightarrow\infty}\|A_{k}-A\|=0$.
\end{theorem}

\section{PROOFS}

The following lemma describes an important property of the
superdiagonal entries in (2).

\begin{lemma}
An upper triangular matrix in $M_{n}(\mathbb{C})$ of the form
$$\left(\begin{array}{ccccc}
\alpha_{_{11}}&\alpha_{_{12}}&\alpha_{_{13}}&\cdots&\alpha_{_{1n}}\\
0&\alpha_{_{22}}&\alpha_{_{23}}&\cdots&\alpha_{_{2n}}\\
0&0&\alpha_{_{33}}&\cdots&\alpha_{_{3n}}\\
\vdots&\vdots&\vdots&\ddots&\vdots\\
0&0&0&\cdots&\alpha_{_{nn}}\\
\end{array}\right)$$ is strongly irreducible if and only if the
equation $\alpha_{_{11}}=\alpha_{_{22}}=\cdots=\alpha_{_{nn}}$ and
the inequality $\alpha_{_{i,i+1}}\neq 0$ for $1\leq i\leq n-1$ both
hold.
\end{lemma}

\begin{proof}
If the matrix is strongly irreducible, then the equation
$\alpha_{_{11}}=\cdots=\alpha_{_{nn}}$ holds. Write $\alpha$ for
$\alpha_{_{ii}}$, $1\leq i\leq n$. Because every strongly
irreducible matrix is similar to a Jordan matrix, we know that there
is an invertible matrix in $M_{n}(\mathbb{C})$ such that
$$\left(\begin{array}{ccccc}
\alpha_{_{11}}&\alpha_{_{12}}&\alpha_{_{13}}&\cdots&\alpha_{_{1n}}\\
0&\alpha_{_{22}}&\alpha_{_{23}}&\cdots&\alpha_{_{2n}}\\
0&0&\alpha_{_{33}}&\cdots&\alpha_{_{3n}}\\
\vdots&\vdots&\vdots&\ddots&\vdots\\
0&0&0&\cdots&\alpha_{_{nn}}\\
\end{array}\right)
\left(\begin{array}{ccccc}
x_{_{11}}&x_{_{12}}&x_{_{13}}&\cdots&x_{_{1n}}\\
x_{_{21}}&x_{_{22}}&x_{_{23}}&\cdots&x_{_{2n}}\\
x_{_{31}}&x_{_{32}}&x_{_{33}}&\cdots&x_{_{3n}}\\
\vdots&\vdots&\vdots&\ddots&\vdots\\
x_{_{n1}}&x_{_{n2}}&x_{_{n3}}&\cdots&x_{_{nn}}\\
\end{array}\right)$$
$$=\left(\begin{array}{ccccc}
x_{_{11}}&x_{_{12}}&x_{_{13}}&\cdots&x_{_{1n}}\\
x_{_{21}}&x_{_{22}}&x_{_{23}}&\cdots&x_{_{2n}}\\
x_{_{31}}&x_{_{32}}&x_{_{33}}&\cdots&x_{_{3n}}\\
\vdots&\vdots&\vdots&\ddots&\vdots\\
x_{_{n1}}&x_{_{n2}}&x_{_{n3}}&\cdots&x_{_{nn}}\\
\end{array}\right)
\left(\begin{array}{ccccc}
\alpha&1&0&\cdots&0\\
0&\alpha&1&\cdots&0\\
0&0&\alpha&\cdots&0\\
\vdots&\vdots&\vdots&\ddots&\vdots\\
0&0&0&\cdots&\alpha\\
\end{array}\right).$$ This equation yields that  $x_{_{ij}}=0$ for
$i>j$ and $x_{_{ii}}=\alpha_{_{i,i+1}}x_{_{i+1,i+1}}$ for $i=1, 2,
\ldots, n-1$. Hence, we obtain
$x_{_{kk}}=\prod^{n-1}_{i=k}\alpha_{_{i,i+1}}x_{_{nn}}$. If
$\alpha_{_{i,i+1}}=0$ for some $i$ in $\{1,2,\ldots,n-1\}$, then the
matrix $(x_{_{ij}})_{_{1\leq i,j\leq n}}$ is not invertible.
Therefore the inequality $\alpha_{_{i,i+1}}\neq 0$ holds for $1\leq
i\leq n-1$.

On the other hand, if $\alpha_{i,i+1}\neq 0$ holds for $i=1, 2,
\ldots, n$, then every matrix in $M_{n}(\mathbb{C})$ commuting with
the matrix $(\alpha_{_{ij}})_{_{1\leq i,j\leq n}}$ can be expressed
in the form
$$X=\left(\begin{array}{ccccc}
x_{_{11}}&x_{_{12}}&x_{_{13}}&\cdots&x_{_{1n}}\\
0&x_{_{11}}&x_{_{23}}&\cdots&x_{_{2n}}\\
0&0&x_{_{11}}&\cdots&x_{_{3n}}\\
\vdots&\vdots&\vdots&\ddots&\vdots\\
0&0&0&\cdots&x_{_{11}}\\
\end{array}\right).$$ If $X$ is an idempotent, then it must be $I$ or
$0$. Thus the matrix $(\alpha_{_{ij}})_{_{1\leq i,j\leq n}}$ is
strongly irreducible.
\end{proof}

Applying this lemma, we obtain the following corollary.

\begin{corollary}
In (2), the function $\phi^{n}_{i,i+1}$ in $L^{\infty}(\mu_{n})$
satisfies $\phi^{n}_{i,i+1}(\lambda)\neq 0$ almost everywhere on
$\Lambda_{n}$ for $i=1,2,\ldots, n-1$.
\end{corollary}

In this corollary, the Multiplication operator
$M_{\phi^{n}_{i,i+1}}$ induced by the function $\phi^{n}_{i,i+1}$ is
not invertible in general. But $M_{\phi^{n}_{i,i+1}}$ can be
approximated by a sequence of invertible Multiplication operators in
$\mathscr{L}(L^{2}(\mu_{n}))$. Meanwhile, replacing the
superdiagonal entries with invertible ones enable us to simplify the
problem. That is why we add the hypothesis (ii) in Theorem 1.1.
Precisely, we obtain the following two lemmas.

\begin{lemma}
If an operator $A_{n}$ is a direct integral of strongly irreducible
operators stated as in (2) in the form
$$A_{n}=\left(\begin{array}{cccccc}
M_{\phi_{n}}&M_{\phi^{n}_{12}}&M_{\phi^{n}_{13}}&\cdots&M_{\phi^{n}_{1n}}\\
0&M_{\phi_{n}}&M_{\phi^{n}_{23}}&\cdots&M_{\phi^{n}_{2n}}\\
0&0&M_{\phi_{n}}&\cdots&M_{\phi^{n}_{3n}}\\
\vdots&\vdots&\vdots&\ddots&\vdots\\
0&0&0&\cdots&M_{\phi_{n}}\\
\end{array}\right)_{n\times n},$$ then for every positive integer
$k$, there exists an operator $A_{nk}$ in the form
$$A_{nk}=\left(\begin{array}{cccccc}
M_{\phi_{n}}&M_{\phi^{n}_{12,k}}&M_{\phi^{n}_{13}}&\cdots&M_{\phi^{n}_{1n}}\\
0&M_{\phi_{n}}&M_{\phi^{n}_{23,k}}&\cdots&M_{\phi^{n}_{2n}}\\
0&0&M_{\phi_{n}}&\cdots&M_{\phi^{n}_{3n}}\\
\vdots&\vdots&\vdots&\ddots&\vdots\\
0&0&0&\cdots&M_{\phi_{n}}\\
\end{array}\right)_{n\times n}$$ with invertible $M_{\phi^{n}_{i,i+1,k}}$
for $1\leq i\leq n-1$ such that $\|A_{n}-A_{nk}\|<\dfrac{1}{k}$.
\end{lemma}

\begin{proof}
For $\lambda$ in $\Lambda_{n}$, we construct
${{\phi^{n}_{i,i+1,k}}}$ in the form
$$\phi^{n}_{{i,i+1,k}}(\lambda)=\left\{\begin{array}{ll}
\phi^{n}_{{i,i+1}}(\lambda),&
\mbox{if}~|\phi^{n}_{{i,i+1}}(\lambda)|\geq\dfrac{1}{kn};\\
\dfrac{\phi^{n}_{{i,i+1}}(\lambda)}{kn|\phi^{n}_{{i,i+1}}(\lambda)|},&
\mbox{if}~0<|\phi^{n}_{{i,i+1}}(\lambda)|<\dfrac{1}{kn};\\
\dfrac{1}{kn},&\mbox{if}~\phi^{n}_{{i,i+1}}(\lambda)=0.\\
\end{array}\right.$$
Thus
$\|M_{{\phi^{n}_{i,i+1}}}-M{_{\phi^{n}_{i,i+1,k}}}\|<\dfrac{1}{k(n-1)},\quad
1\leq i\leq n-1.$ Therefore we obtain
$$\|A_{n}-A_{nk}\|\leq\sum^{n-1}_{i=1}
\|M_{{\phi^{n}_{i,i+1}}}-M{_{\phi^{n}_{i,i+1,k}}}\|
<(n-1)\dfrac{1}{k(n-1)}=\dfrac{1}{k}.$$ By the definition, the
operator $M_{\phi^{n}_{i,i+1,k}}$ is invertible for $1\leq i\leq
n-1$.
\end{proof}

\begin{lemma}
If an operator $A_{n}$ is a direct integral of strongly irreducible
operators stated as in (2) in the form
$$A_{n}=\left(\begin{array}{cccccc}
M_{\phi_{n}}&M_{\phi^{n}_{12}}&M_{\phi^{n}_{13}}&\cdots&M_{\phi^{n}_{1n}}\\
0&M_{\phi_{n}}&M_{\phi^{n}_{23}}&\cdots&M_{\phi^{n}_{2n}}\\
0&0&M_{\phi_{n}}&\cdots&M_{\phi^{n}_{3n}}\\
\vdots&\vdots&\vdots&\ddots&\vdots\\
0&0&0&\cdots&M_{\phi_{n}}\\
\end{array}\right)_{n\times n}$$ such that $M_{\phi^{n}_{i,i+1}}$
is invertible in $\mathscr{L}(L^{2}(\mu_{n}))$ for $i=1,2,\ldots,
n-1$, then there exists an invertible operator $X_{n}$ in
$\mathscr{L}((L^{2}(\mu_{n}))^{(n)})$ such that
$X_{n}A_{n}X^{-1}_{n}$ is in the form
$$X_{n}A_{n}X^{-1}_{n}=\left(\begin{array}{cccccc}
M_{\phi_{n}}&I&0&\cdots&0\\
0&M_{\phi_{n}}&I&\cdots&0\\
0&0&M_{\phi_{n}}&\cdots&0\\
\vdots&\vdots&\vdots&\ddots&\vdots\\
0&0&0&\cdots&M_{\phi_{n}}\\
\end{array}\right). \eqno(3)$$
\end{lemma}

\begin{proof}
We construct an invertible upper triangular operator-valued matrix
$X_{n}$ in $\mathscr{L}((L^{2}(\mu_{n}))^{(n)})$ as follows.

Choose an invertible operator $M_{f^{}_{nn}}$ in
$\mathscr{L}(L^{2}(\mu_{n}))$. Fix an operator $M_{f^{}_{ii}}$ by
the equation
$$M_{\phi^{n}_{i,i+1}}M_{f^{}_{i+1,i+1}}=M_{f^{}_{ii}}$$ for
$i=1,\ldots,n-1$. Let $\{M_{f^{}_{ii}}\}^{n}_{i=1}$ be the main
diagonal ($0$-diagonal) entries of $X_{n}$. Notice that every
operator in the set $\{M_{f^{}_{ii}}\}^{n}_{i=1}$ is invertible in
$\mathscr{L}(L^{2}(\mu_{n}))$.

Choose an operator $M_{f^{n}_{n-1,n}}$ in
$\mathscr{L}(L^{2}(\mu_{n}))$. Fix an operator $M_{f^{n}_{i,i+1}}$
by the equation
$$M_{\phi^{n}_{i,i+1}}M_{f^{n}_{i+1,i+2}}+
M_{\phi^{n}_{i,i+2}}M_{f^{n}_{i+2,i+2}}=M_{f^{n}_{i,i+1}}$$ for
$i=1,\ldots,n-2$. Let $\{M_{f^{n}_{i,i+1}}\}^{n-1}_{i=1}$ be the
$1$-diagonal entries of $X_{n}$.

Choose an operator $M_{f^{n}_{n-l,n}}$ in
$\mathscr{L}(L^{2}(\mu_{n}))$, where $l$ is a positive integer such
that $1\leq l\leq n-1$. Fix an operator $M_{f^{n}_{i,i+l}}$ by the
equation
$$M_{\phi^{n}_{i,i+1}}M_{f^{n}_{i+1,i+l+1}}+
M_{\phi^{n}_{i,i+2}}M_{f^{n}_{i+2,i+l+1}}+\cdots+
M_{\phi^{n}_{i,i+l+1}}M_{f^{n}_{i+l+1,i+l+1}}=M_{f^{n}_{i,i+l}}$$
for $i=1,\ldots,n-l-1$. Let $\{M_{f^{n}_{i,i+l}}\}^{n-l}_{i=1}$ be
the $l$-diagonal entries of $X_{n}$.

Choose an operator $M_{f^{n}_{1n}}$ in $\mathscr{L}(L^{2}(\mu_{n}))$
to be the $n$-diagonal entry of $X_{n}$.

Therefore we obtain an invertible operator-valued matrix $X_{n}$ in
the form
$$X_{n}=\left(\begin{array}{cccccc}
M_{f_{11}}&M_{f^{n}_{12}}&M_{f^{n}_{13}}&\cdots&M_{f^{n}_{1n}}\\
0&M_{f_{22}}&M_{f^{n}_{23}}&\cdots&M_{f^{n}_{2n}}\\
0&0&M_{f_{33}}&\cdots&M_{f^{n}_{3n}}\\
\vdots&\vdots&\vdots&\ddots&\vdots\\
0&0&0&\cdots&M_{f_{nn}}\\
\end{array}\right)$$ such that the equation (3) holds.
\end{proof}

By Lemma 2.4, we can reduce equation (2) to the form
$$A=
\bigoplus\limits^{\infty}_{n=1}\int\limits_{\Lambda_n}
\left(\begin{array}{cccccc}
M_{\phi_{n}}&I&0&\cdots&0\\
0&M_{\phi_{n}}&I&\cdots&0\\
0&0&M_{\phi_{n}}&\cdots&0\\
\vdots&\vdots&\vdots&\ddots&\vdots\\
0&0&0&\cdots&M_{\phi_{n}}\\
\end{array}\right)_{n\times n}(\lambda) d\mu_{n}(\lambda),\eqno (4)$$
in the sense of similar equivalence.

For a regular Borel measure $\nu$ on $\mathbb{C}$ with compact
support $K$, define $N_{\nu}$ on $L^{2}(\nu)$ by $N_{\nu}f=z\cdot f$
for each $f$ in $L^{2}(\nu)$.

\begin{lemma}
Let an operator $A_{n}$ be in the form
$$A_{n}=\left(\begin{array}{cccccc}
N^{(\infty)}_{\nu_{n}}&I&0&\cdots&0\\
0&N^{(\infty)}_{\nu_{n}}&I&\cdots&0\\
0&0&N^{(\infty)}_{\nu_{n}}&\cdots&0\\
\vdots&\vdots&\vdots&\ddots&\vdots\\
0&0&0&\cdots&N^{(\infty)}_{\nu_{n}}\\
\end{array}\right)_{n\times n},$$ where $\nu_{n}$ is
a regular Borel measure and supported on some compact set $K_{n}$
such that $0<\nu_{n}(K_{n})<\infty$. Then the strongly irreducible
decomposition of $A_{n}$ is not unique up to similarity.
\end{lemma}

\begin{proof}
To prove this lemma, we need to construct two bounded maximal
abelian sets of idempotents in $\{A_{n}\}^{\prime}$ which are not
similar.

We can write $N^{(\infty)}_{\nu_{n}}$ in the form
$N_{\nu_{n}}\otimes I_{l^{2}}$, where $I_{l^{2}}$ is the identity
operator on $l^{2}$. Denote by $\mathscr{P}$ the set of all the
spectral projections of $N_{\nu_{_n}}$. This set forms a bounded
maximal abelian set of idempotents in $\{N_{\nu_{_n}}\}^{\prime}$.
Let $\{e_{k}\}^{\infty}_{k=1}$ be an orthonormal basis for $l^{2}$.
Denote by $E_{k}$ the projection such that
$\textrm{ran}E_{k}=\{\lambda e_{k}:\lambda\in\mathbb{C}\}$. Let
$\mathscr {Q}_{1}\triangleq\{P\in\mathscr{L}(l^{2}):P=P^{*}=P^{2}
\in\{E_{k}:k\in\mathbb{N}\}^{\prime\prime}\}$. Denote by $\chi_{_S}$
the characteristic function for a Borel subset $S$ in the interval
$[0,1]$ and let $\hat{\mathscr
{Q}}_{2}\triangleq\{M_{\chi_{_S}}\in\mathscr{L}(L^{2}[0,1]):
S\subset[0,1]\mbox{ is Borel.}\}$. There is a unitary operator
$U:L^{2}[0,1]\rightarrow l^{2}$ such that
$UPU^{*}\in\mathscr{L}(l^{2})$ for every $P\in\hat{\mathscr
{Q}}_{2}$. The sets ${\mathscr {Q}}_{2}\triangleq U\hat{\mathscr
{Q}}_{2}U^{*}$ and $\mathscr {Q}_{1}$ are two bounded maximal
abelian sets of idempotents in $\mathscr{L}(l^{2})$ but they are not
unitarily equivalent. By the fact that $W^{*}(\mathscr {P})\otimes
W^{*}(\mathscr {Q}_{1})$ and $W^{*}(\mathscr {P})\otimes
W^{*}(\mathscr {Q}_{2})$ are both maximal abelian von Neumann
algebras, we obtain that $$\mathscr {F}_{1}\triangleq\{P\in
W^{*}(\mathscr {P})\otimes W^{*}(\mathscr {Q}_{1}):P=P^{*}=P^{2}\}$$
and $$\mathscr {F}_{2}\triangleq\{P\in W^{*}(\mathscr {P})\otimes
W^{*}(\mathscr {Q}_{2}):P=P^{*}=P^{2}\}$$ are both maximal abelian
sets of idempotents in $\{N_{\nu_{n}}\otimes
I_{l^{2}}\}^{\prime}=L^{\infty}(\nu_{n})\otimes\mathscr{L}(l^{2})$.

We need to prove that $\mathscr {F}^{(n)}_{i}$ is a bounded maximal
abelian set of idempotents in $\{A_{n}\}^{\prime}$ for $i=1,2$.

An operator $X$ in $\{A_{n}\}^{\prime}$ can be expressed in the form
$$X=\left(\begin{array}{cccccc}
X_{11}&X_{12}&X_{13}&\cdots&X_{1n}\\
X_{21}&X_{22}&X_{23}&\cdots&X_{2n}\\
X_{31}&X_{32}&X_{33}&\cdots&X_{3n}\\
\vdots&\vdots&\vdots&\ddots&\vdots\\
X_{n1}&X_{n2}&X_{n3}&\cdots&X_{nn}\\
\end{array}\right)_{n\times n}. \eqno (5)$$ We prove that $X_{ij}$
is in $\{N_{\nu_{n}}\otimes I_{l^{2}}\}^{\prime}$. Note that
$\mathscr{P}^{(\infty)}$ is the set of all the spectral projections
of $N_{\nu_{n}}\otimes I_{l^{2}}$. Fix an projection $P$ in
$(\mathscr{P}^{(\infty)})^{(n)}$. The operator $A$ can be expressed
in the form $$A_{n}=A_{n1}\oplus A_{n2},$$ where
$$A_{ni}=\left(\begin{array}{cccccc}
N^{(\infty)}_{\nu_{{ni}}}&I&0&\cdots&0\\
0&N^{(\infty)}_{\nu_{{ni}}}&I&\cdots&0\\
0&0&N^{(\infty)}_{\nu_{{ni}}}&\cdots&0\\
\vdots&\vdots&\vdots&\ddots&\vdots\\
0&0&0&\cdots&N^{(\infty)}_{\nu_{{ni}}}\\
\end{array}\right),\ i=1,2.$$ The measures $\nu_{{n1}}$ and
$\nu_{{n2}}$ are mutually singular and their supports depend on the
characteristic functions corresponding to $P$ and $I-P$. Hence $X$
can also be expressed in the form
$$X=\begin{pmatrix}
Y_{11}&Y_{12}\\
Y_{21}&Y_{22}\\
\end{pmatrix}
\begin{array}{l}
{\textrm{ran}}P\\
{\textrm{ran}}(I-P)\\
\end{array}.$$
The equations $A_{n1}Y_{12}=Y_{12}A_{n2}$ and
$A_{n2}Y_{21}=Y_{21}A_{n1}$ yield that $Y_{12}=Y_{21}=0$. Therefore
$P$ reduces $X$ and $X_{ij}$s are in $\{N_{\nu_{n}}\otimes
I_{l^{2}}\}^{\prime}$. A computation shows that the equation
$X_{ij}=0$ holds for $i>j$ and $X_{ii}=X_{11}$ for $i=2,\ldots,n$ in
(5). Furthermore, if $X$ as in (5) is an idempotent, then so is
every main diagonal entry $X_{ii}$ of $X$.

We assume that $X$ is an idempotent in $\{A_{n}\}^{\prime}$ and
commutes with $\mathscr {F}^{(n)}_{1}$. Hence $X_{ii}$ commutes with
$\mathscr {F}_{1}$. The fact that $\mathscr {F}_{1}$ is a maximal
abelian set of idempotents implies that $X_{ii}$ belongs to
$\mathscr {F}_{1}$. Thus $X_{ii}$ commutes with $X_{ij}$. For the
$1$-diagonal entries, the equation $2X_{ii}X_{i,i+1}-X_{i,i+1}=0$
yields $X_{i,i+1}=0$, for $i=1,\ldots,n-1$. By this way, the
$k$-diagonal entries of $X$ are all zero, for $k=2,\ldots,n$.
Therefore $X$ is in $\mathscr {F}^{(n)}_{1}$. Both $\mathscr
{F}^{(n)}_{1}$ and $\mathscr {F}^{(n)}_{2}$ are bounded maximal
abelian sets of idempotents in $\{A_{n}\}^{\prime}$.

We prove that $\mathscr {F}^{(n)}_{1}$ and $\mathscr {F}^{(n)}_{2}$
are not similar in $\{A_{n}\}^{\prime}$. Every operator $X$ in
$\{A_{n}\}^{\prime}$ can be written in the form
$$X=\int_{\sigma(N_{\nu_{_n}})} X(\lambda) d\nu_{n}(\lambda).$$

Suppose that there is an invertible operator $X$ in
$\{A_{n}\}^{\prime}$ such that $$X\mathscr
{F}^{(n)}_{2}X^{-1}=\mathscr {F}^{(n)}_{1}.$$ For each $P$ in
$\mathscr {F}^{(n)}_{2}$, the projection $P(\lambda)$ is either of
rank $\infty$ or $0$, for almost every $\lambda$ in
$\sigma(N_{\nu_{n}})$. But there exists an projection $Q$ in
$\mathscr {F}^{(n)}_{1}$ such that $Q(\lambda)$ is of rank $n$, for
almost every $\lambda$ in $\sigma(N_{\nu_{n}})$. This is a
contradiction. Therefore $\mathscr {F}^{(n)}_{1}$ and $\mathscr
{F}^{(n)}_{2}$ are not similar in $\{A_{n}\}^{\prime}$.
\end{proof}

By (\cite{Shi_2}, Theorem 3.3), we have the following corollary.

\begin{corollary}
Let an operator $A_{n}$ be in the form
$$A_{n}=\left(\begin{array}{cccccc}
N^{(m)}_{\nu_{n}}&I&0&\cdots&0\\
0&N^{(m)}_{\nu_{n}}&I&\cdots&0\\
0&0&N^{(m)}_{\nu_{n}}&\cdots&0\\
\vdots&\vdots&\vdots&\ddots&\vdots\\
0&0&0&\cdots&N^{(m)}_{\nu_{n}}\\
\end{array}\right)_{n\times n},$$ where $m$ is a positive integer
and $\nu_{n}$ is a regular Borel measure supported on some compact
set $K_{n}$ such that $0<\nu_{n}(K_{n})<\infty$. Then the strongly
irreducible decomposition of $A_{n}$ is unique up to similarity.
\end{corollary}

For a regular Borel measure $\nu$ with compact support, Denote by
$J_{m}(\nu)$ an operator in the form
$$J_{m}(\nu)=\left(\begin{array}{cccccc}
N_{\nu}&I&0&\cdots&0\\
0&N_{\nu}&I&\cdots&0\\
0&0&N_{\nu}&\cdots&0\\
\vdots&\vdots&\vdots&\ddots&\vdots\\
0&0&0&\cdots&N_{\nu}\\
\end{array}\right)_{m\times m}
\begin{array}{c}
L^{2}(\nu)\\
L^{2}(\nu)\\
L^{2}(\nu)\\
\vdots\\
L^{2}(\nu)\\
\end{array}. \eqno (6)$$

\begin{lemma}
Every operator in $\{J_{m}(\nu)\}^{\prime}$ is in the form
$$\begin{pmatrix}
M_{\phi_{_{1}}}&M_{\phi_{_{2}}}&\cdots&M_{\phi_{_{m}}}\\
0&M_{\phi_{_{1}}}&\cdots&M_{\phi_{_{m-1}}}\\
\vdots&\vdots&\ddots&\vdots\\
0&0&\cdots&M_{\phi_{_{1}}}\\
\end{pmatrix},$$
where $\phi_{{i}}$ is in $L^{\infty}(\nu)$ for $1\leq i\leq m$.
\end{lemma}

\begin{proof}
By a similar computation as in Lemma 2.5, we obtain that every
operator in $\{J_{m}(\nu)\}^{\prime}$ is in the form
$$\begin{pmatrix}
M_{\phi_{_{11}}}&M_{\phi_{_{12}}}&\cdots&M_{\phi_{_{1m}}}\\
0&M_{\phi_{_{22}}}&\cdots&M_{\phi_{_{2,m-1}}}\\
\vdots&\vdots&\ddots&\vdots\\
0&0&\cdots&M_{\phi_{_{mm}}}\\
\end{pmatrix}.$$
By the equation
$$N_{\nu}M_{\phi_{_{i,j+1}}}+M_{\phi_{_{i+1,j+1}}}
=M_{\phi_{_{i,j}}}+M_{\phi_{_{i,j+1}}}N_{\nu},$$ the $k$-diagonal
entries are as required for $1\leq k\leq n$.
\end{proof}

\begin{lemma}
Let $m_{1}$ and $m_{2}$ be two positive integers such that
$m_{1}>m_{2}$. Then the following equations hold:
\begin{enumerate}
\item
$\{B\in\mathscr{L}((L^{2}(\nu))^{(m_{2})},(L^{2}(\nu))^{(m_{1})}):
J_{m_{_1}}(\nu)B=BJ_{m_{_2}}(\nu)\}\\=
\{(C^{T},0)^{T}:C\in\{J_{m_{_2}}(\nu)\}^{\prime}\}$.
\item
$\{B\in\mathscr{L}((L^{2}(\nu))^{(m_{1})},(L^{2}(\nu))^{(m_{2})}):
J_{m_{_2}}(\nu)B=BJ_{m_{_1}}(\nu)\}\\=
\{(0,C):C\in\{J_{m_{_2}}(\nu)\}^{\prime}\}$.
\end{enumerate}
\end{lemma}

\begin{proof}
We only need to prove the first equation. The second equation can be
obtained by the same method. Let $B=(C^{T},D^{T})^{T}$ such that $C$
and $D$ are in the form
$$C=\begin{pmatrix}
B_{11}&B_{12}&\cdots&B_{1m_{_2}}\\
B_{21}&B_{22}&\cdots&B_{2m_{_2}}\\
\vdots&\vdots&\ddots&\vdots\\
B_{m_{_2}1}&B_{m_{_2}2}&\cdots&B_{m_{_2}m_{_2}}\\
\end{pmatrix},\quad
D=\begin{pmatrix}
B_{m_{_2}+1,1}&\cdots&B_{m_{_2}+1,m_{_2}}\\
\vdots&&\vdots\\
B_{m_{_1}1}&\cdots&B_{m_{_1}m_{_2}}\\
\end{pmatrix}.$$
By a similar computation as in Lemma 2.5, we can obtain that
$P^{(m_{_1})}B=BP^{(m_{_2})}$ for every spectral projection $P$ of
$N_{\nu}$. Thus, every $B_{ij}$ belongs to $\{N_{\nu}\}^{\prime}$,
for $1\leq i\leq m_{1}$ and $1\leq j\leq m_{2}$. For
$i=1,\ldots,m_{1}-1$, the equation
$N_{\nu}B_{i1}+B_{i+1,1}=B_{i1}N_{\nu}$ yields $B_{i+1,1}=0$. For
$i=2,\ldots,m_{1}-1$, the equation
$N_{\nu}B_{i2}+B_{i+1,2}=B_{i2}N_{\nu}$ yields $B_{i+1,2}=0$. By
this way, we can obtain $B_{ij}=0$ for $i>j$. Hence $D=0$ and a
further computation shows that $C\in\{J_{m_{_2}}(\nu)\}^{\prime}$.
\end{proof}

\begin{lemma}
Let $m_{1}$, $m_{2}$, $r_{1}$, $r_{2}$ be positive integers. If an
idempotent $P$ in $M_{n}(\mathbb{C})$ is in the form
$$P=\begin{pmatrix}
I_{m_{_1}}&0&R_{11}&R_{12}\\
0&0_{r_{_1}}&R_{21}&R_{22}\\
0&0&I_{m_{_2}}&0\\
0&0&0&0_{r_{_2}}\\
\end{pmatrix},$$
where $I_{m}$ is the identity operator in $M_{m}(\mathbb{C})$, then
there exists an invertible operator $X$ in the form
$$X=\begin{pmatrix}
I_{m_{_1}}&0&0&R_{12}\\
0&I_{r_{_1}}&-R_{21}&0\\
0&0&I_{m_{_2}}&0\\
0&0&0&I_{r_{_2}}\\
\end{pmatrix}~~ \mbox{and}~~
X^{-1}=\begin{pmatrix}
I_{m_{_1}}&0&0&-R_{12}\\
0&I_{r_{_1}}&R_{21}&0\\
0&0&I_{m_{_2}}&0\\
0&0&0&I_{r_{_2}}\\
\end{pmatrix}$$
such that
$$XPX^{-1}=\begin{pmatrix}
I_{m_{_1}}&0&0&0\\
0&0_{r_{_1}}&0&0\\
0&0&I_{m_{_2}}&0\\
0&0&0&0_{r_{_2}}\\
\end{pmatrix}.$$
\end{lemma}

Note that the equation $P^{2}=P$ implies that $R_{11}=R_{22}=0$ and
the construction of $X$ depends on $P$. In the following example, we
construct an operator $A$ and prove the strongly irreducible
decomposition of $A$ is unique up to similarity.

\begin{example}
Let $A=J^{(2)}_{3}(\nu)\oplus J^{(3)}_{2}(\nu)\oplus
J^{(2)}_{1}(\nu)$. We prove that for every two bounded maximal
abelian sets of idempotents $\mathscr{P}$ and $\mathscr{Q}$ in
$\{A\}^{\prime}$, there is an invertible operator $X$ in
$\{A\}^{\prime}$ such that the equation
$\mathscr{P}=X\mathscr{Q}X^{-1}$ holds and
$$V(\{A\}^{\prime})=\{f\mbox{ is bounded
Borel}:\sigma(N_{\nu})\rightarrow
\mathbb{N}\oplus\mathbb{N}\oplus\mathbb{N}\},$$
$$K_{0}(\{A\}^{\prime})=\{f\mbox{ is bounded
Borel}:\sigma(N_{\nu})\rightarrow
\mathbb{Z}\oplus\mathbb{Z}\oplus\mathbb{Z}\}.$$

Denote by $\mathscr{P}_{m_{_1}}$ the set of all the idempotents in
$\{J_{m_{_1}}(\nu)\}^{\prime}$. Note that $\mathscr{P}_{m_{_1}}$
equals the set of all the spectral projections of
$N^{(m_{_1})}_{\nu}$. Denote by $\mathscr{E}_{m_{_2}}$ the set of
all the diagonal projections in $M_{m_{_2}}(\mathbb{C})$ and by
$\mathscr{F}_{m_{_1},m_{_2}}$ the set of all the projections in
$\{\mathscr{P}_{m_{_1}}\otimes\mathscr{E}_{m_{_2}}\}^{\prime\prime}$.

Let $\mathscr{P}=
\mathscr{F}_{3,2}\oplus\mathscr{F}_{2,3}\oplus\mathscr{F}_{1,2}$. We
can verify that $\mathscr{P}$ is a bounded maximal abelian set of
idempotents in $\{A\}^{\prime}$. Then we only need to prove that for
every bounded maximal abelian set of idempotents $\mathscr{Q}$ in
$\{A\}^{\prime}$, there is an invertible operator $X$ in
$\{A\}^{\prime}$ such that $\mathscr{P}=X\mathscr{Q}X^{-1}$.

We reduce the rest into two claims:
\begin{enumerate}
\item For every idempotent $P$ in $\{A\}^{\prime}$,
there is an invertible operator $X$ in $\{A\}^{\prime}$ such that
$XPX^{-1}$ belongs to $\mathscr{P}$.
\item There are seven idempotents $\{Q_{k}\}^{7}_{k=1}$ in
$\mathscr{Q}$ such that for almost every $\lambda$ in
$\sigma(N_{\nu})$, $\{Q_{k}(\lambda)\}^{7}_{k=1}$ and
$\mathscr{Q}(\lambda)$ generate the same bounded maximal abelian set
of idempotents.
\end{enumerate}
Every operator $B$ in $\{A\}^{\prime}$ can be expressed in the form
$$B=\begin{pmatrix}
B^{11}&B^{12}&\cdots&B^{17}\\
B^{21}&B^{22}&\cdots&B^{27}\\
\vdots&\vdots&\ddots&\vdots\\
B^{71}&B^{72}&\cdots&B^{77}\\
\end{pmatrix},$$ where
$$\begin{array}{lll}
B^{11}=\begin{pmatrix}
b^{11}_{1}&b^{11}_{2}&b^{11}_{3}\\
0&b^{11}_{1}&b^{11}_{2}\\
0&0&b^{11}_{1}\\
\end{pmatrix},&
B^{13}=\begin{pmatrix}
b^{13}_{1}&b^{13}_{2}\\
0&b^{13}_{1}\\
0&0\\
\end{pmatrix},&
B^{16}=\begin{pmatrix}
b^{16}_{1}\\
0\\
0\\
\end{pmatrix},\\
\\
B^{31}=\begin{pmatrix}
0&b^{31}_{1}&b^{31}_{2}\\
0&0&b^{31}_{1}\\
\end{pmatrix},&
B^{33}=\begin{pmatrix}
b^{33}_{1}&b^{33}_{2}\\
0&b^{33}_{1}\\
\end{pmatrix},&
B^{36}=\begin{pmatrix}
b^{36}_{1}\\
0\\
\end{pmatrix},\\
\\
B^{61}=\begin{pmatrix}
0&0&b^{31}_{1}\\
\end{pmatrix},&
B^{63}=\begin{pmatrix}
0&b^{33}_{1}\\
\end{pmatrix},&
B^{66}=\begin{pmatrix}
b^{66}_{1}\\
\end{pmatrix},
\end{array}$$
other $B_{ij}$s are expressed as follows:
\begin{itemize}
\item For $1\leq i,j\leq 2$, $B^{ij}$s are of the same form;
\item For $3\leq i,j\leq 5$, $B^{ij}$s are of the same form;
\item For $6\leq i,j\leq 7$, $B^{ij}$s are of the same form;
\item For $1\leq i,j\leq 2$ and $3\leq j\leq 5$, $B^{ij}$s are of
the same form;
\item For $1\leq i,j\leq 2$ and $6\leq j\leq 7$, $B^{ij}$s are of
the same form;
\item For $3\leq i,j\leq 5$ and $1\leq j\leq 2$, $B^{ij}$s are of
the same form;
\item For $3\leq i,j\leq 5$ and $6\leq j\leq 7$, $B^{ij}$s are of
the same form;
\item For $6\leq i,j\leq 7$ and $1\leq j\leq 2$, $B^{ij}$s are of
the same form;
\item For $6\leq i,j\leq 7$ and $3\leq j\leq 5$, $B^{ij}$s are of
the same form,
\end{itemize}
where $b^{ij}_{k}$s belong to $\{N_{\nu}\}^{\prime}$ for $1\leq
i,j\leq 7$ and $1\leq k\leq 3$.

For $B$ expressed in the above form, there is a unitary operator
$U_{1}$ such that
$$U_{1}BU^{*}_{1}=\begin{pmatrix}
B_{11}&B_{12}&B_{13}\\
0&B_{22}&B_{23}\\
0&0&B_{33}\\
\end{pmatrix},$$ where $B_{ij}$s are in the form
$$B_{11}=\begin{pmatrix}
b^{11}_{1}&b^{12}_{1}&\vdots&b^{13}_{1}&b^{14}_{1}&b^{15}_{1}&\vdots&b^{16}_{1}&b^{17}_{1}\\
b^{21}_{1}&b^{22}_{1}&\vdots&b^{23}_{1}&b^{24}_{1}&b^{25}_{1}&\vdots&b^{26}_{1}&b^{27}_{1}\\
\cdots&\cdots&&\cdots&\cdots&\cdots&&\cdots&\cdots\\
0&0&\vdots&b^{33}_{1}&b^{34}_{1}&b^{35}_{1}&\vdots&b^{36}_{1}&b^{37}_{1}\\
0&0&\vdots&b^{43}_{1}&b^{44}_{1}&b^{45}_{1}&\vdots&b^{46}_{1}&b^{47}_{1}\\
0&0&\vdots&b^{53}_{1}&b^{54}_{1}&b^{55}_{1}&\vdots&b^{56}_{1}&b^{57}_{1}\\
\cdots&\cdots&&\cdots&\cdots&\cdots&&\cdots&\cdots\\
0&0&\vdots&0&0&0&\vdots&b^{66}_{1}&b^{67}_{1}\\
0&0&\vdots&0&0&0&\vdots&b^{76}_{1}&b^{77}_{1}\\
\end{pmatrix}_{7\times 7},\eqno (7)$$
$$B_{12}=\begin{pmatrix}
b^{11}_{2}&b^{12}_{2}&\vdots&b^{13}_{2}&b^{14}_{2}&b^{15}_{2}\\
b^{21}_{2}&b^{22}_{2}&\vdots&b^{23}_{2}&b^{24}_{2}&b^{25}_{2}\\
\cdots&\cdots&&\cdots&\cdots&\cdots\\
0&0&\vdots&b^{33}_{2}&b^{34}_{2}&b^{35}_{2}\\
0&0&\vdots&b^{43}_{2}&b^{44}_{2}&b^{45}_{2}\\
0&0&\vdots&b^{53}_{2}&b^{54}_{2}&b^{55}_{2}\\
\cdots&\cdots&&\cdots&\cdots&\cdots\\
0&0&\vdots&0&0&0\\
0&0&\vdots&0&0&0\\
\end{pmatrix}_{7\times 5},\quad
B_{13}=\begin{pmatrix}
b^{11}_{3}&b^{12}_{3}\\
b^{21}_{3}&b^{22}_{3}\\
\cdots&\cdots\\
0&0\\
0&0\\
0&0\\
\cdots&\cdots\\
0&0\\
0&0\\
\end{pmatrix}_{7\times 2},$$
$$B_{22}=\begin{pmatrix}
b^{11}_{1}&b^{12}_{1}&\vdots&b^{13}_{1}&b^{14}_{1}&b^{15}_{1}\\
b^{21}_{1}&b^{22}_{1}&\vdots&b^{23}_{1}&b^{24}_{1}&b^{25}_{1}\\
\cdots&\cdots&&\cdots&\cdots&\cdots\\
0&0&\vdots&b^{33}_{1}&b^{34}_{1}&b^{35}_{1}\\
0&0&\vdots&b^{43}_{1}&b^{44}_{1}&b^{45}_{1}\\
0&0&\vdots&b^{53}_{1}&b^{54}_{1}&b^{55}_{1}\\
\end{pmatrix},\quad
B_{23}=\begin{pmatrix}
b^{11}_{2}&b^{12}_{2}\\
b^{21}_{2}&b^{22}_{2}\\
\cdots&\cdots\\
0&0\\
0&0\\
0&0\\
\end{pmatrix},$$
$$B_{33}=\begin{pmatrix}
b^{11}_{1}&b^{12}_{1}\\
b^{21}_{1}&b^{22}_{1}\\
\end{pmatrix}.$$

If $P$ is an idempotent in $\{U_{1}AU^{*}_{1}\}^{\prime}$, then by
the proof of (\cite{Shi_2}, Lemma 3.4), we can construct an
invertible operator $X$ in $\{U_{1}AU^{*}_{1}\}^{\prime}$ of the
form
$$X=\begin{pmatrix}
X_{1}&&\\
&X_{2}&\\
&&X_{3}\\
\end{pmatrix}\oplus
\begin{pmatrix}
X_{1}&\\
&X_{2}\\
\end{pmatrix}\oplus
X_{1},$$ such that
$$XPX^{-1}=\begin{pmatrix}
P_{11}&P_{12}&P_{13}\\
0&P_{22}&P_{23}\\
0&0&P_{33}\\
\end{pmatrix},$$ where the main diagonal blocks as in (7) are
diagonal projections. There is also a unitary operator $U_{2}$ in
$\{U_{1}AU^{*}_{1}\}^{\prime}$ of the form
$$U_{2}=\begin{pmatrix}
U^{2}_{1}&&\\
&U^{2}_{2}&\\
&&U^{2}_{3}\\
\end{pmatrix}\oplus
\begin{pmatrix}
U^{2}_{1}&\\
&U^{2}_{2}\\
\end{pmatrix}\oplus
U^{2}_{1},$$ such that the equation
$$\begin{pmatrix}
U^{2}_{1}&&\\
&U^{2}_{2}&\\
&&U^{2}_{3}\\
\end{pmatrix}P_{11}\begin{pmatrix}
U^{2}_{1}&&\\
&U^{2}_{2}&\\
&&U^{2}_{3}\\
\end{pmatrix}^{*}(\lambda)$$
$$=\begin{pmatrix}
I_{s_{_1}}&0&\vdots&*&*&\vdots&*&*\\
0&0_{t_{_1}}&\vdots&*&*&\vdots&*&*\\
\cdots&\cdots&&\cdots&\cdots&&\cdots&\cdots\\
0&0&\vdots&I_{s_{_2}}&0&\vdots&*&*\\
0&0&\vdots&0&0_{t_{_2}}&\vdots&*&*\\
\cdots&\cdots&&\cdots&\cdots&&\cdots&\cdots\\
0&0&\vdots&0&0&\vdots&I_{s_{_3}}&0\\
0&0&\vdots&0&0&\vdots&0&0_{t_{_3}}\\
\end{pmatrix}$$ holds for almost every $\lambda$ in
$\sigma(N_{\nu})$, where $s_{i}$ and $t_{i}$ are non-negative
integers.

Write $U_{2}XPX^{-1}U^{*}_{2}$ in the form
$$\widehat{P}=U_{2}XPX^{-1}U^{*}_{2}=\begin{pmatrix}
\widehat{P}_{11}&\widehat{P}_{12}&\widehat{P}_{13}\\
0&\widehat{P}_{22}&\widehat{P}_{23}\\
0&0&\widehat{P}_{33}\\
\end{pmatrix}.$$ By Lemma 2.9, we can construct an invertible
operator $Y_{1}$ in $\{U_{1}AU^{*}_{1}\}^{\prime}$ such that
$\widehat{P}_{11}$, $\widehat{P}_{22}$, and $\widehat{P}_{33}$
become diagonal projections after similar transformation.
Furthermore, we can construct an invertible operator $Y_{2}$ in
$\{U_{1}AU^{*}_{1}\}^{\prime}$ such that the $1$-diagonal blocks of
$Y_{1}\widehat{P}Y^{-1}_{1}$ vanish after similar transformation.
And then we can construct an invertible operator $Y_{3}$ in
$\{U_{1}AU^{*}_{1}\}^{\prime}$ such that the $2$-diagonal blocks of
$Y_{2}Y_{1}\widehat{P}Y^{-1}_{1}Y^{-1}_{2}$  vanish after similar
transformation. Thus we finish the proof of claim (i).

To prove claim (ii), we need to define a $\nu$-measurable function
$r_{_Q}$ with respect to an idempotent $Q$ in $\{A\}^{\prime}$.
Without loss of generality, we assume that $Q\sim (P_{31}\oplus
P_{32})\oplus(P_{21}\oplus P_{22}\oplus P_{23})\oplus (P_{11}\oplus
P_{12})\in\mathscr{F}_{3,2}\oplus\mathscr{F}_{2,3}\oplus\mathscr{F}_{1,2}$.
Define
$$\begin{array}{ccl}
r_{_Q}(\lambda)&\triangleq&\dfrac{1}{3}\textrm{Tr}(P_{31}(\lambda)+P_{32}(\lambda))\\
&&+\dfrac{1}{2}\textrm{Tr}(P_{21}(\lambda)+P_{22}(\lambda)+P_{22}(\lambda))\\
&&+\ \ \
\textrm{Tr}(P_{11}(\lambda)+P_{12}(\lambda)),~\lambda\in{\sigma}(N_{\nu}),
\end{array}$$
where Tr stands for the standard trace of a square matrix. Note that
$r_{_Q}$ stays invariant up to similarity. By the proof of
(\cite{Shi_2}, Lemma 3.5), we can obtain that there are seven
idempotents $Q_{i}$ in $\mathscr{Q}$ such that
\begin{itemize}
\item the equation $r_{Q_{i}}(\lambda)=1$ holds a.~e.\ on $\sigma(N_{\nu})$
for $i=1, \ldots, 7$ and
\item the equation $Q_{i}Q_{j}=0$ holds for $i\neq j$.
\end{itemize}
The idempotent $Q_{i}$ may be in the form
$$Q_{i}(\lambda)\sim\left\{\begin{array}{ll}
I_{3}\oplus 0,&\lambda\in\Lambda_{3};\\
I_{2}\oplus 0,&\lambda\in\Lambda_{2};\\
I_{1}\oplus 0,&\lambda\in\Lambda_{1},\\
\end{array}\right.$$ where $\{\Lambda_{i}\}^{3}_{i=1}$ is a Borel
partition of $\sigma(N_{\nu})$. We can choose finitely many spectral
projections of $N^{(7)}_{\nu}$ to cut $Q_{i}$s and to piece together
new $Q_{i}$s such that every $Q_{i}$ belongs to a $\mathscr{P}_{m}$
for $m=1,2,3$. We finish the proof of claim (ii).

By the proof of (\cite{Shi_2}, Lemma 3.6) and the idempotents
$\{Q_{i}\}^{7}_{i=1}$ constructed above, we can obtain an invertible
operator $X$ in $\{A\}^{\prime}$ such that
$X\mathscr{Q}X^{-1}=\mathscr{P}$. Therefore the strongly irreducible
decomposition of $A$ is unique up to similarity.

Assume that $Q_{1}$, $Q_{2}$ and $Q_{3}$ are in $\mathscr{P}_{1}$,
$\mathscr{P}_{2}$ and $\mathscr{P}_{3}$ respectively. Then there is
a group isomorphism $\alpha$ such that $\alpha([Q_{1}])=(1,0,0)$,
$\alpha([Q_{2}])=(0,1,0)$, and $\alpha([Q_{3}])=(0,0,1)$, where
$[Q_{i}]$ stands for the similar equivalence class of $Q_{i}$ in
$V(\{A\}^{\prime})=\bigcup^{\infty}_{n=1}M_{n}(\{A\}^{\prime})\slash\sim$.
Thus we obtain $\alpha([I])=(2,3,2)$, where $I$ is the identity
operator in $\{A\}^{\prime}$. Furthermore, a routine computation
yields that $V(\{A\}^{\prime})$ and $K_{0}(\{A\}^{\prime})$ are of
the forms at the beginning of this example.

\begin{proof}[Proof of Theorem 1.1]

By the calculating in the above example, we can prove that the
strongly irreducible decomposition of
$\bigoplus^{k}_{i=1}J^{(n_{i})}_{m_{_i}}(\nu)$ is unique up to
similarity, where $ m_{_i}$, $n_{_i}$ and $k$ are all positive
integers.

There is a unitary operator $V$ such that $VAV^{*}$ can be expressed
in the form as described at the beginning of Example 2.10. Then we
can apply the above lemmas to perform calculation as we need. Note
that the equation
$$\{(\bigoplus^{k_{1}}_{i=1}J^{(n_{i})}_{m_{_i}}(\nu_{1}))
\oplus(\bigoplus^{k_{2}}_{j=1}J^{(n_{j})}_{m_{_j}}(\nu_{2}))\}^{\prime}
=\{\bigoplus^{k_{1}}_{i=1}J^{(n_{i})}_{m_{_i}}(\nu_{1})\}^{\prime}
\oplus\{\bigoplus^{k_{2}}_{j=1}J^{(n_{j})}_{m_{_j}}(\nu_{2})\}^{\prime}$$
holds for mutually singular Borel measures $\nu_{1}$ and $\nu_{2}$.

By Lemma 2.5, if the strongly irreducible decomposition of $A$ is
unique up to similarity, then every multiplicity function
$m_{_{{\phi_{_n}}}}$ is bounded. Then we can obtain that
$V(\{A\}^{\prime})$ and $K_{0}(\{A\}^{\prime})$ are as described in
the theorem. On the other hand, if the strongly irreducible
decomposition of $A$ is not unique up to similarity, then there is a
number $m$ in $\{n\in\mathbb{N}:\mu(\Lambda_{n})>0\}$ such that the
multiplicity function $m_{_{{\phi_{_m}}}}$ takes $\infty$ in its
codomain on a Borel subset $\Gamma_{m1}$ of measure nonzero in its
domain. Therefore in $K_{0}(\{A\}^{\prime})$, every Borel function
$f$ vanishes on $\Gamma_{m1}$. This is a contradiction.
\end{proof}

The proof of Theorem 1.2 is an application of Lemma 2.3.

\end{example}


\begin{acknowledgements}
The author is grateful to Professor Chunlan Jiang and Professor
Guihua Gong for their advice and comments on writing this paper.
Also the author was supported in part by NSFC Grant (No. 10731020)
and NSFC Grant (No.10901046).
\end{acknowledgements}


\end{document}